\documentclass[a4paper,11pt]{article}
\usepackage{amssymb}
\usepackage{amsmath}
\usepackage{latexsym,a4wide}
\input xypic
\input xy
\xyoption{all}

\newtheorem{example}{Example}[section]
\newtheorem{definition}[example]{Definition}
\newtheorem{proposition}[example]{Proposition}
\newtheorem{prop}[example]{Proposition}

\newtheorem{thm}[example]{Theorem}

\newenvironment{proof}{\noindent \textbf{Proof:} }{$\Box$ \mbox{}}

\newcommand{\llabto}[2]{\stackrel{#2}{\rule[0.52ex]{#1em}{0.099ex}\hspace{-0.4em}\longrightarrow}}

\begin{document}

\date{}
\author{ \"{O}zg\"{u}n  G\"{u}rmen Alansal and Erdal Ulualan}
\title{PEIFFER ELEMENTS IN THE MOORE COMPLEX OF A BISIMPLICIAL GROUP}
\maketitle

\begin{abstract}
In this work, we explain `Peiffer pairings'  in the Moore (bi)complex of a \textit{bisimplicial group} and give their applications for crossed modules and crossed squares.
\\


\end{abstract}

\section*{Introduction}
 In \cite{mutpor2}, Mutlu and Porter generalized Peiffer elements to higher dimensions giving systematic maps  of generating them. They used ideas based on the works of Conduch{\'e} \cite{con} and Carrasco-Cegarra \cite{car} and gave images of the functions $F_{\alpha,\beta}$ in the Moore complex of a simplicial group. For further accounts  of  these functions in simplicial  Lie and commutative algebras see  \cite{AA} and \cite{Arvasi, Arvasi2}, respectively. The purpose of this paper  is to give their results for  \textit{bisimplicial groups}.

Crossed modules were introduced by Whitehead in \cite{W}. Crossed $n$-cubes (cf. \cite{els}) model
connected $(n+1)$-types and for $ n=2, $ they are just considered as crossed
squares (cf. \cite{lg}). These crossed squares are related to 2-crossed modules introduced by Conduch\'e in \cite{con}.

In a letter to Brown and Loday, dated in the mid 1980s, Conduch\'e pointed out that the mapping cone complex of a crossed square, constructed by Loday in \cite{loday}, has a 2-crossed module structure. As a generalization of this result, Conduch\'e has proven Theorem 0.1 in \cite{cond}. In this work, he showed that crossed squares are related to bisimplicial groups in the same way crossed modules are related to simplicial groups. In this paper we aim  to shed some light on the 2-crossed modules and crossed squares structures given by Conduch\'e, by use of the functions $F_{\alpha,\beta}$  in the Moore (bi)complex of a bisimplicial group. In particular, we see  in Proposition \ref{3} that the $h$-map of the crossed square given by Conduch\'e corresponds to the functions $F_{\alpha,\beta}$. There is a crossed module version  of this  in Proposition \ref{2}. These results, by a similar way, can be iterated to multisimplicial (or $n$-simplicial) groups and also the relationship between crossed $n$-cubes  and $n$-simplicial groups can be constructed in terms of these functions. In this case, a 1-truncated $n$-simplicial group  gives   a crossed $n$-cube structure.

\section{(Bi)simplicial groups}
In this section, the definitions and properties related to (bi)simplicial groups, some of which are classical, are recalled. We refer the reader to May's and Loday's books (cf. \cite{may} and
\cite{loday2}) and Artin-Mazur's \cite{Artin} article  for most of the basic
properties of (bi)simplicial groups  that we will be needing.

Let $\Delta$ be the category of finite ordinals $[n]=\{0<1<\cdots
<n\}$. A \textit{simplicial group} is a functor from the opposite
category $\Delta ^{op}$ to the category of groups $\mathbf{Grp}$.
That is, a simplicial group $\mathbf{G}$ consists of a family of
groups $G_{n}$ together with homomorphisms $d_i^n:G_n\rightarrow
G_{n-1}$, $0\leqslant i \leqslant n$, $(n\neq 0)$ and $s_j^n:G_n
\rightarrow G_{n+1}$, $0\leqslant j \leqslant n$, called face and
degeneracy maps, satisfying the usual simplicial identities given in \cite{dusk, may}.

Given a simplicial group $\mathbf{G}$, \textit{the Moore complex }($\mathbf{%
NG}$,$\partial $) of $\mathbf{G}$ is the normal chain complex
defined by
\begin{equation*}
(NG)_{n}=\ker d_{0}^{n}\cap\ker d_{1}^{n}\cap \dots\cap \ker
d_{n-1}^{n}
\end{equation*}%
with boundaries $\partial _{n}:NG_{n}\rightarrow NG_{n-1}$ induced
from $d_{n}^{n}$ by restriction.

A $2$-simplicial group or \textit{a bisimplicial group} $\mathbf{G_{\ast,\ast}}$ is a functor from the
product category $\Delta ^{op}\times \Delta ^{op}$ to the category
of groups {\bf Grp}, with the face and degeneracy maps given by
\begin{align*}
d_{i}^{h^{(pq)}}:G_{p,q}\rightarrow G_{p-1,q} &  \\
s_{i}^{h^{(pq)}}:G_{p,q}\rightarrow G_{p+1,q} & \quad 0\leq i\leq p \\
d_{j}^{v^{(pq)}}:G_{p,q}\rightarrow G_{p,q-1} &  \\
s_{j}^{v^{(pq)}}:G_{p,q}\rightarrow G_{p,q+1} & \quad 0\leq j\leq q
\end{align*}
such that the maps $d_{i}^{h^{(pq)}},s_{i}^{h^{(pq)}}$ commute with $d_{j}^{v^{(pq)}},s_{j}^{v^{(pq)}}$
and that $d_{i}^{h^{(pq)}},s_{i}^{h^{(pq)}}$ (resp. $d_{j}^{v^{(pq)}},s_{j}^{v^{(pq)}}$) satisfy the
usual simplicial identities.

We think of $d_{j}^{v^{(pq)}},s_{j}^{v^{(pq)}}$ as the vertical operators and $%
d_{i}^{h^{(pq)}},s_{i}^{h^{(pq)}}$ as the horizontal operators. If $\mathbf{G_{\ast,\ast}}$
is a bisimplicial group, it is convenient to think of an element of
$G_{p,q}$ as a product of a $p$-simplex and a $q$-simplex.

 The Moore bicomplex of a bisimplicial group $\mathbf{G_{\ast,\ast}}$ is defined  by
$$
NG_{n,m}=\underset{(i,j)=(0,0)}{
\overset{(n-1,m-1)}{\bigcap }}\text{Ker}
d_{i}^{h^{(nm)}}\cap \text{Ker}d_{j}^{v^{(nm)}}
$$
with the boundary homomorphisms
$$
\partial _{i}^{h^{(nm)}}:NG_{n,m}\longrightarrow
NG_{n-1,m}
$$
and
$$
\partial_{j}^{v^{(nm)}}:NG_{n,m}\longrightarrow
NG_{n,m-1}
$$
induced by the face maps $d_{i}^{h^{(nm)}}$ and $d_j^{v^{(nm)}}$ where $0\leqslant i\leqslant n \neq 0$, $0\leqslant j\leqslant m \neq 0$.

This Moore bicomplex is illustrated by the following diagram.
$$
\xymatrix{&\ar[d]&\ar[d]&\ar[d]\\
\cdots \ar[r]& NG_{2,2} \ar[d] \ar[r]^{\vdots}&NG_{1,2}\ar[r]_{\partial_1^{h^{(12)}}}^{\vdots}\ar[d]_{\partial_2^{v^{(12)}}}&NG_{0,2}\ar[d]^{\partial_2^{v^{(02)}}}
\\ \cdots \ar[r]&NG_{2,1}\ar[d] \ar[r]&NG_{1,1}\ar[d]_{\partial_1^{v^{(11)}}}\ar[r]_{\partial_1^{h^{(11)}}}&NG_{0,1}\ar[d]^{\partial_1^{v^{(01)}}}\\
\cdots \ar[r]& NG_{2,0}\ar[r]_{\partial_2^{h^{(20)}}}&NG_{1,0}\ar[r]_{\partial_1^{h^{(10)}}}&NG_{0,0}.}
$$
This definition extends easily to multisimplicial groups.

\section{Peiffer pairings for bisimplicial groups}

This section is devoted to description  of the functions $F_{\alpha,\beta}$ in the Moore complex of a bisimplicial group.
The following notation and terminology is derived from \cite{car}. For the ordered set $[n]=\{0<1<\dots <n\}$, let $\alpha
_{i}^{n}:[n+1]\rightarrow \lbrack n]$ be the increasing surjective
map given by;
\begin{equation*}
\alpha _{i}^{n}(j)=\left\{
\begin{array}{ll}
j & \text{if }j\leqslant
 i, \\
j-1 & \text{if }j>i.%
\end{array}%
\right.
\end{equation*}
Let $S(n,n-r)$ be the set of all monotone increasing surjective maps from $%
[n]$ to $[n-r]$. This can be generated from the various $\alpha
_{i}^{n}$ by composition. The composition of these generating maps
is subject to the following rule: $\alpha _{j}\alpha _{i}=\alpha
_{i-1}\alpha _{j},j<i$. This
implies that every element $\alpha \in S(n,n-r)$ has a unique expression as $%
\alpha =\alpha _{i_{1}}\circ \alpha _{i_{2}}\circ \dots \circ \alpha
_{i_{r}}$ with $0\leqslant
 i_{1}<i_{2}< \dots <i_{r}\leqslant
 n-1$, where the
indices $i_{k}$ are the elements of $[n]$ such that
$\{i_{1},\dots,i_{r}\}=\{i:\alpha (i)=\alpha
(i+1)\}$. We thus can identify $S(n,n-r)$ with the set $%
\{(i_{r},\dots,i_{1}):0\leqslant
 i_{1}<i_{2}<\dots<i_{r}\leqslant
 n-1\}$. In
particular, the single element of $S(n,n)$, defined by the identity
map on $[n]$, corresponds to the empty 0-tuple ( ) denoted by
$\emptyset _{n}$. Similarly the only element of $S(n,0)$ is
$(n-1,n-2,\dots,0)$. For all $n\geqslant  0$, let
\begin{equation*}
S(n)=\bigcup_{0\leqslant r\leqslant n}S(n,n-r).
\end{equation*}
For example
$$S(0)=\{\emptyset\}, S(1)=\{\emptyset, (0)\}, S(2)=\{\emptyset, (1),(0),(1,0)\}.$$

The following  terminology is derived from \cite{cond}. For $n,q\in {\Bbb N}$ with $q\leqslant n$, let $S(n,q)$ be the set of nondecreasing surjections from
$[n]$ to $[q]$ as defined above. For $\sigma\in S(n,q)$ the target of $\sigma $ is called $b(\sigma):q=b(\sigma)$. The set $S(n)$ is partially ordered set by the following relation: $\sigma\leqslant\tau $ if, for $i\in[n]$, one has $\sigma(i)\geqslant \tau(i)$ where $[b(\sigma)]$ and $[b(\tau)]$ are considered as subsets of ${\Bbb N}$.

Given $r\in {\Bbb N}^*$ and $\underline{n}=(n_1,\dots, n_r)\in {\Bbb N}^r$, let $S(\underline{n})=S(n_1)\times \cdots \times S(n_r)$ with the product partial order.
For $\sigma=(\sigma_1,\dots,\sigma_r)\in S(\underline{n})$ let $b(\sigma)=(b(\sigma_1),\dots, b(\sigma_r))$.

The following result was proved by Conduch\'e in \cite{condd}.
\begin{thm}
Let $G_*$ be a $r$-simplicial group. Then, for $\underline{n}\in {\Bbb N}^r $ , the group $G_{\underline{n}}$ is a $S(\underline{n})$-semi-direct product of the family of subgroups $(s_\sigma(N(G)_{b(\sigma)}))_{\sigma\in S(\underline{n})}$, where $s_{\sigma}$ is the composite of the iterated degeneracies $s_{\sigma_i}$ corresponding to each simplicial structure.
\end{thm}

Now we can give the functions $F_{\alpha,\beta}$ in the Moore complex of a bisimplicial group.

Given $\underline{n}=\left( k_{1},k_{2}\right) \in {\Bbb
N}\times{\Bbb N}$, let $S\left( \underline{n}\right) =S\left(
k_{1}\right) \times S\left( k_{2}\right) $ with the product
(partial) order. Let $\underline{\alpha },\ \underline{\beta }\in
S\left( \underline{n}\right) $
and $\underline{\alpha }=\left( \alpha _{1},\alpha _{2}\right) ;\ \underline{%
\beta }=\left( \beta _{1},\beta _{2}\right) $ where $ \alpha _{i},
\beta_{i} \in S\left(k_{i}\right) $ for $1\leqslant i \leqslant 2.$  The pairings that
we will need
\[
\left\{ F_{\underline{\alpha },\underline{\beta }}:NG_{\underline{n}-\#%
\underline{\alpha }}\times NG_{\underline{n}-\#\underline{\beta }%
}\longrightarrow NG_{\underline{n}}\ \ ;\ \underline{\alpha },\underline{%
\beta }\in S\left( \underline{n}\right), \underline{\alpha }\neq\underline{\beta } \right\}
\]
are given by composing of the maps in the following diagram
$$
\xymatrix{  &&NG_{k_1- \# \alpha_1 , k_2- \# \alpha_2 }\times
NG_{k_1- \# \beta_1 ,k_2- \# \beta_2}
\ar[dd]_{(s_{\underline{\alpha}},s_{\underline{\beta}})} \ar[rr]^-{F_{\underline{\alpha },\underline{\beta }}} &&
NG_{k_1,k_2}\\ \\
  &&G_{k_1,k_2}\times G_{k_1,k_2} \ar[rr]_{\mu} && G_{k_1,k_2}\ar[uu]_{p}  &&    }
$$
where $s_{\underline{\alpha}}:s_{\alpha _{1}}^{h}s_{\alpha
_{2}}^{v},\ $ and where $s_{\alpha_1}^h=s_{i_{r}}^h\ldots s_{i_{1}}^h$ for $
\alpha_{1}=(i_r,\ldots,i_1)\in S(k_1) $ and similarly $s_{\underline{\beta}}=s_{\beta _{1}}^{h}s_{\beta
_{2}}^{v},\ $  $s_{\beta_2}^v=s_{j_{m}}^v\ldots s_{j_{1}}^v$ for $
\beta_{1}=(j_m,\ldots,j_1)\in S(k_2) $ and  where
$$
p:G_{k_1,k_2} \rightarrow  NG_{k_1,k_2}
$$
 is defined by the composite projection%
\[
p=\left( p_{k_{1}-1}^{h}...p_{0}^{h}\right) \left(
p_{k_{2}-1}^{v}...p_{0}^{v}\right)
\]
where $p_{j}\left( x\right)
=x^{-1}s_{j}d_{j}\left( x\right)$ in both horizontal and vertical directions
and $\mu $ is given by the commutator map.

Thus for  $ \underline{\alpha}=(\alpha_1,\alpha_2),
\underline{\beta}=(\beta_1,\beta_2)\in S(k_1)\times S(k_2) $, we
obtain
\begin{eqnarray*}
F_{\underline{\alpha},\underline{\beta}}(x,y)&=&p\mu
(s_{\underline{\alpha}},s_{\underline{\beta}})(x,y)\\
&=&p\mu
(s^{h}_{\alpha_{1}}s^{v}_{{\alpha_{2}}}(x),s^{h}_{{\beta_{1}}}s^{v}_{{{\beta}_{2}}}(y))\\
&=&p[s^{h}_{{{\alpha}_{1}}}s^{v}_{{{\alpha}_{2}}}(x),s^{h}_{{{\beta}_{1}}}s^{v}_{{{\beta}_{2}}}(y)]
\end{eqnarray*}
where $ x\in NG_{k_1- \# \alpha_1 , k_2- \# \alpha_2 } $ and $ y\in
NG_{k_1- \# \beta_1 ,k_2- \# \beta_2}.$

The construction given above can be easily extended to the $n$-simplicial groups as follows.

Given $n\neq 0, n\in {\Bbb N}$ and $\underline{n}=\left(
k_{1},k_{2},\dots ,k_{n}\right) \in {\Bbb N}\times {\Bbb N}\times
...\times {\Bbb N}={\Bbb N}^n$, let $S\left( \underline{n}\right)
=S\left( k_{1}\right) \times S\left( k_{2}\right) \times ...\times
S\left( k_{n}\right) $ with the product (partial) order as given in \cite{cond}.

Let $\underline{\alpha },\underline{\beta }\in S\left(
\underline{n}\right) $ and $\underline{\alpha }=\left( \alpha
_{1},\alpha _{2},\dots ,\alpha _{n}\right) ;\underline{\beta }=\left(
\beta _{1},\beta _{2},\dots ,\beta _{n}\right) $ where $\alpha _{i}\in
S\left( k_{i}\right) $ and $\beta _{j}\in S\left( k_{j}\right)
,~1\leqslant i,j\leqslant n.$

 The pairings that
we will need
\[
\left\{ F_{\underline{\alpha },\underline{\beta }}:NG_{\underline{n}-\#%
\underline{\alpha }}\times NG_{\underline{n}-\#\underline{\beta }%
}\longrightarrow NG_{\underline{n}}\ \ ;\ \underline{\alpha },\underline{%
\beta }\in S\left( \underline{n}\right) \right\}
\]
are given as composites by the diagram
$$ \xymatrix{  &&NG_{k_{1}- \# \alpha_1 , k_{2}- \# \alpha_2 ,\dots , k_{n}- \#
\alpha_n
 }\times NG_{k_{1}- \# \beta_1 ,k_{2}- \# \beta_2 ,\dots ,k_{n}- \# \beta_n}
\ar[dd]_{(s^{t_1}_{\alpha_1}s^{t_2}_{\alpha_2}...s^{t_n}_{\alpha_n},s^{t_1}_{\beta_1}s^{t_2}_{\beta_2}...s^{t_n}_{\beta_n})}
\ar[rr]^{\ \ \ \ \ \ \ \ \ \ \ \ \ \ \ \ \ \ \ \ \ \ \ \ \ \
F_{\alpha,\beta}} &&
NG_{k_{1},k_{2},\dots ,k_{n}}\\ \\
  &&G_{k_{1},k_{2},\dots ,k_{n}}\times G_{k_{1},k_{2},\dots ,k_{n}} \ar[rr]_{\mu} && G_{k_{1},k_{2},\dots ,k_{n}}\ar[uu]_{p}  &&    }
$$
where $s_{\underline{\alpha}}:s_{\alpha _{1}}^{t_{1}}s_{\alpha
_{2}}^{t_{2}}...s_{\alpha _{n}}^{t_{n}},\ $for $1\leqslant i\leqslant n;\
s_{\alpha _{i}}^{t_{i}}:s_{i_r}^{t_{i}}\cdots s_{i_1}^{t_{i}}$ for $\alpha_i=(i_r,\dots, i_1)\in S(k_i)$ and
similarly $s_{\beta _{i}}^{t_{i}},$ and $p$ is defined by the
composite projection
\[
p=\left( p_{k_{1}-1}^{t_{1}}...p_{0}^{t_{1}}\right) \left(
p_{k_{2}-1}^{t_{2}}...p_{0}^{t_{2}}\right) ...\left(
p_{k_{n}-1}^{t_{n}}...p_{0}^{t_{n}}\right)
\]
where $p_{j}^{t_{k}}\left( x\right)
=x^{-1}s_{j}^{t_{k}}d_{j}^{t_{k}}\left( x\right) ,\ $for any $j$ and
$k$, and where each $t_k$ for $1\leqslant k\leqslant n$ indicates the directions of $n$-simplicial group, and $\mu $ is given by the commutator map.

\section{Calculations of the functions $F_{\alpha,\beta}$ in low dimensions} \label{table}

For $ 0\leqslant k_1, k_2 \leqslant 2 $ we consider the sets
$ S(k_1)\times S(k_2). $  We shall calculate the images
of the functions $F_{\underline{\alpha}, \underline{\beta}}$ for all
$\underline{\alpha}, \underline{\beta}\in S(k_1)\times S(k_2).$

First, consider $(n,m)=(0,1)$ or $(n,m)=(1,0)$. Then, we get  the $F_{\underline{\alpha}, \underline{\beta}}$ functions whose codomain $NG_{0,1}$ or $NG_{1,0}$ respectively. First let $(n,m)=(0,1)$. Then
$$S(n,m)=S(0)\times S(1)=\{(\emptyset,\emptyset),(\emptyset,(0))\}.$$
Let $\underline{\alpha}=(\emptyset,\emptyset)$ and $\underline{\beta}=(\emptyset,(0))$. Then the function $$F_{(\emptyset,\emptyset),(\emptyset,(0))}:NG_{0,1}\times NG_{0,0}\longrightarrow NG_{0,1}$$
can be given as follows:
\begin{align*}
F_{(\emptyset,\emptyset),(\emptyset,(0))}(x, y)&=p\mu(s^h_{\emptyset}s^v_{\emptyset}(x), s^h_{\emptyset}s^v_{(0)}(y))\\
&=p_0^v[id(x),s^v_{0}(y)] \ (\because  s^h_{\emptyset}=s^v_{\emptyset}=id)\\
&=[x,s_0^{v^{(00)}}y]s_0^{v^{(00)}}d_0^{v^{(01)}}[s_0^{v^{(00)}}y,x]\\
&=[x,s_0^{v^{(00)}}y][s_0^{v^{(00)}}y,s_0^{v^{(00)}}d_0^{v^{(01)}}(x)]\\
&=[x,s_0^{v^{(00)}}y][s_0^{v^{(00)}}y,1]\ ((\because x\in\ker d_0^{^v{(01)}}=NG_{0,1})\\
&=[x,s_0^{v^{(00)}}y]
\end{align*}
for $x\in NG_{0,1}$ and $y\in NG_{0,0}$.

Suppose now that $(n,m)=(1,0)$. Then we take $S(1)\times S(0)=\{(\emptyset,\emptyset),((0),\emptyset)\}$. Let  $\underline{\alpha}=(\emptyset,\emptyset)$ and $\underline{\beta}=((0),\emptyset)$. Then the function
$$F_{(\emptyset,\emptyset),((0),\emptyset)}:NG_{1,0}\times NG_{0,0}\longrightarrow NG_{1,0}$$
is defined by
\begin{align*}
F_{(\emptyset,\emptyset),((0),\emptyset)}(x, y)&=p\mu(s^h_{\emptyset}s^v_{\emptyset}(x), s^h_{0}s^v_{\emptyset}(y))\\
&=p_0^h[id(x),s^h_{0}(y)]\\
&=[x,s_0^{h^{(00)}}y]s_0^{h^{(00)}}d_0^{h^{(10)}}[s_0^{h^{(00)}}y,x]\\
&=[x,s_0^{h^{(00)}}y][s_0^{h^{(00)}}y,s_0^{h^{(00)}}d_0^{h^{(10)}}(x)]\\
&=[x,s_0^{h^{(00)}}y][s_0^{h^{(00)}}y,1]\ ((\because x\in\ker d_0^{^h{(10)}}=NG_{1,0})\\
&=[x,s_0^{h^{(00)}}y]
\end{align*}
for all $x\in NG_{1,0}$ and $y\in NG_{00}$.

We give the calculations of other functions listed below in Appendix A.

$\hspace{-1cm}
\begin{tabular}{|l|l|c|c|l|}
\hline
 $ $&$F_{\underline{\alpha} ,\underline{\beta} }:\text{Domain} \rightarrow \text{Codomain}$& $\underline{\alpha} $ & $\underline{\beta} $ & $F_{\underline{\alpha} ,\underline{\beta} }\left( x,y\right) $ \\
\hline
\multicolumn{1}{|l|}{$1.$}& {$ NG_{0,1}\times NG_{0,0}\rightarrow NG_{0,1}$} &
$\left( \emptyset,\emptyset \right) $ & $\left( \emptyset ,(0)\right) $ & $[x,s_{0}^{v^{(11)}}(y)]$ \\
\hline
\multicolumn{1}{|l|}{$2.$}& {$NG_{1,0}\times NG_{0,0}\rightarrow NG_{1,0}$} &
$\left( \emptyset,\emptyset\right)$ & $\left( (0),\emptyset \right)$ & $[x,s_0^{v^{(00)}}y]$
 \\
\hline
\multicolumn{1}{|l|}{$3.$} &{$NG_{1,1}\times NG_{1,0}\rightarrow NG_{1,1}$} &
$\left( \emptyset ,\emptyset\right) $ & $\left( \emptyset ,(0)\right) $ &$[x,s_0^{v^{(10)}}y]$ \\
\hline
\multicolumn{1}{|l|}{$4.$} &{$NG_{1,1}\times NG_{0,1}\rightarrow NG_{1,1}$} &
 $\left( \emptyset,\emptyset \right)$ & $\left( (0),\emptyset \right)$ & $[x,s_{0}^{h^{(01)}}(y) ]$ \\
\hline
\multicolumn{1}{|l|}{$5.$} &{$NG_{1,1}\times NG_{0,0}\rightarrow NG_{1,1}$} &
$\left( \emptyset,\emptyset\right)$ & $\left( (0),(0) \right)$ & $[ x ,s_0^{v^{(00)}}s_{0}^{h^{(01)}}( y) ]$ \\
\hline
\multicolumn{1}{|l|} {$6.$ }&{$NG_{0,1}\times NG_{1,0}\rightarrow NG_{1,1}$} &
$\left( (0),\emptyset \right) $ & $\left( \emptyset ,(0)\right) $ & $[s_{0}^{h^{(01)}}x,s_{0}^{v^{(10)}}y]$ \\
\hline
\multicolumn{1}{|l|} {$7.$} &{$NG_{0,1}\times NG_{0,1}\rightarrow NG_{0,2}$} &
$\left( \emptyset, (0) \right) $ & $\left( \emptyset,(1) \right) $ &$[ s_{0}^{v^{(01)}}x,s_{1}^{v^{(01)}}y][s_{1}^{v^{(01)}}y,s_{1}^{v^{(01)}}x] $ \\
\hline
\multicolumn{1}{|l|}{$8.$} &{$NG_{1,0}\times NG_{1,0}\rightarrow NG_{2,0}$} &
$( (0),\emptyset)$ & $\left( (1),\emptyset \right)$ & $[s_{0}^{h^{(10)}}x,s_{1}^{h^{(10)}}y][s_{1}^{h^{(10)}}y,s_{1}^{h^{(10)}}x]$ \\
  \hline
\multicolumn{1}{|l|}{$9.$} &{$NG_{1,1}\times NG_{1,1}\rightarrow NG_{1,2}$} &
$\left( \emptyset ,(0)\right)$ & $\left( \emptyset,(1) \right)$ & $[s_{0}^{v^{(11)}}x,s_{1}^{v^{(11)}}y][s_{1}^{v^{(11)}}y,s_{1}^{v^{(11)}}x]$ \\
\hline
\multicolumn{1}{|l|}{$10.$} &{$NG_{1,1}\times NG_{0,2}\rightarrow NG_{1,2}$} &
$\left( \emptyset, (1) \right) $ & $\left( (0),\emptyset\right) $ & $[s_{1}^{v^{(11)}}x,s_{0}^{h^{(02)}}y]$ \\
\hline
\multicolumn{1}{|l|}{$11.$} &{$NG_{1,1}\times NG_{0,2}\rightarrow NG_{1,2}$} &
$\left( \emptyset, (0) \right) $ & $\left( (0),\emptyset\right) $ & $[s_{0}^{v^{(11)}}x,s_{0}^{h^{(02)}}y]$ \\
\hline
\multicolumn{1}{|l|}{$12.$} &{$NG_{0,1}\times NG_{1,1}\rightarrow NG_{1,2}$} &
$\left(  (0),(1) \right) $ & $\left(\emptyset, (0)\right) $ & $[s_{0}^{h^{(02)}}s_{1}^{v^{(01)}}x,s_{0}^{v^{(11)}}y]$ \\
\hline
\multicolumn{1}{|l|}{$13.$} &{$NG_{0,1}\times NG_{1,1}\rightarrow NG_{1,2}$} &
$\left(  (0),(0) \right) $ & $\left(\emptyset, (1)\right) $ & $[s_{0}^{h^{(02)}}s_{0}^{v^{(01)}}x,s_{1}^{v^{(11)}}y]$ \\
\hline
\multicolumn{1}{|l|} {$14.$} &{$NG_{1,1}\times NG_{1,1}\rightarrow NG_{2,1}$} &
$\left( (0),\emptyset \right) $ & $\left((1),\emptyset\right) $ & $[s_{0}^{h^{(11)}}x,s_{1}^{h^{(11)}}y][s_{1}^{h^{(11)}}y,s_{1}^{h^{(11)}}x]$ \\
\hline
\multicolumn{1}{|l|} {$15.$} &{$NG_{1,1}\times NG_{2,0}\rightarrow NG_{2,1}$} &
$\left(  (1),\emptyset \right) $ & $\left(\emptyset,(0)\right) $ & $[s_{1}^{h^{(11)}}x,s_{0}^{v^{(20)}}y]$ \\
\hline
\multicolumn{1}{|l|}{$16.$} &{$NG_{1,1}\times NG_{2,0}\rightarrow NG_{2,1}$} &
$\left(  (0),\emptyset \right) $ & $\left(\emptyset,(0)\right) $ & $[s_{0}^{h^{(11)}}x,s_{0}^{v^{(20)}}y]$ \\
\hline
\multicolumn{1}{|l|}{$17.$}&{$NG_{1,0}\times NG_{1,1}\rightarrow NG_{2,1}$} &
$\left( (1), (0) \right) $ & $\left((0),\emptyset\right) $ & $[s_1^{h^{(11)}}s_{0}^{v^{(10)}}x,s_{0}^{h^{(11)}}y]$ \\
\hline
\multicolumn{1}{|l|}{$18.$} &{$NG_{1,0}\times NG_{1,1}\rightarrow NG_{2,1}$} &
$\left( (0), (0) \right) $ & $\left((1),\emptyset\right) $ & $[s_0^{h^{(11)}}s_{0}^{v^{(10)}}x,s_{1}^{h^{(11)}}y]$ \\
\hline
\multicolumn{1}{|l|} {$19.$} &{$NG_{1,2}\times NG_{1,2}\rightarrow NG_{2,2}$} &
$\left( (0),\emptyset \right) $ & $\left((1),\emptyset\right) $ & $[s_{0}^{h^{(12)}}x,s_{1}^{h^{(12)}}y][s_{1}^{h^{(12)}}y,s_{1}^{h^{(12)}}x]$ \\
\hline
\multicolumn{1}{|l|} {$20.$} &{$NG_{2,1}\times NG_{2,1}\rightarrow NG_{2,2}$} &
$(\emptyset ,(0)) $ & $(\emptyset ,(1)) $ & $[s_{0}^{v^{(21)}}(x),s_{1}^{v^{(21)}}(y)][s_{1}^{v^{(21)}}(y),s_{1}^{v^{(21)}}(x)]$ \\
\hline
\multicolumn{1}{|l|}{$21.$} &{$NG_{1,2}\times NG_{2,1}\rightarrow NG_{2,2}$} &
$\left(  (1),\emptyset \right) $ & $\left(\emptyset,(0)\right) $ & $[s_{1}^{h^{(12)}}(x),s_{0}^{v^{(21)}}(y)]$ \\
\hline
\multicolumn{1}{|l|}{$22.$} &{$NG_{1,2}\times NG_{2,1}\rightarrow NG_{2,2}$} &
$\left(  (0),\emptyset \right) $ & $\left(\emptyset,(0)\right) $ & $[s_{0}^{h^{(12)}}(x),s_{0}^{v^{(21)}}(y)]$ \\
\hline
\multicolumn{1}{|l|}{$23.$}&{$NG_{1,1}\times NG_{2,1}\rightarrow NG_{2,2}$} &
$\left( (0), (1) \right) $ & $\left(\emptyset,(0)\right) $ & $[s_{0}^{v^{(21)}}s_{1}^{h^{(11)}}(x),s_{0}^{v^{(21)}}(y)]
$ \\
\hline
\multicolumn{1}{|l|}{$24.$}&{$NG_{1,1}\times NG_{2,1}\rightarrow NG_{2,2}$} &
$\left( (0), (0) \right) $ & $\left(\emptyset,(1)\right) $ & $[s_{1}^{v^{(21)}}s_{0}^{h^{(11)}}(x),s_{1}^{v^{(21)}}(y)]$ \\
\hline
\end{tabular}%
$
\\

\begin{definition}
Let $\mathbf{G}_{\ast ,\ast }$ be a bisimplicial group and  $n,m>1,$ and $%
D_{n,m}$ the subgroup in $G_{n,m}$ generated by degenerate elements. Let $%
N_{n,m}$ be the normal subgroup of $G_{n,m}$ generated by elements
of the form
$$F_{\underline{\alpha },\underline{\beta }}(x,y)\ \text{with}\ \ \underline{\alpha}=(\alpha_1,\alpha_2),
\underline{\beta}=(\beta_1,\beta_2)\in S(n)\times S(m)$$ where $
x\in NG_{n- \# \alpha_1 , m- \# \alpha_2 } $ and $ y\in NG_{n-
\# \beta_1 ,m- \# \beta_2}.$
\end{definition}

Mutlu and Porter in  \cite{mutpor2} have defined a normal subgroup  $N_n$ of $NG_n$
generated by elements of the forms $F_{\alpha, \beta}(x,y).$
Furthermore they proved that there is an
equality
\begin{equation*}
\partial _{n}(NG_{n}\cap D_{n})=\partial _{n}(N_{n}^{G}\cap D_{n}).
\end{equation*}
By a  similar way  we can write the following equalities
$$
\partial_n ^h (NG_{n,m}\cap D_{n,m})=\partial_n ^h (N_{n,m}\cap D_{n,m}),
$$
and
$$
\partial_m ^v (NG_{n,m}\cap D_{n,m})=\partial_m^v (N_{n,m}\cap D_{n,m})
$$
in each direction.

Thus we have the following result.
\begin{prop}
Let $\mathbf{G}_{\ast ,\ast }$ be a bisimplicial group. Then for
$n\geqslant 1, m\geqslant 2$ and $I, J\subseteq [m-1] $ with $I\cup J=[m-1]$, there is the inclusion
\begin{equation*}
\left[ K_{I}\cap K_{H},K_{J}\cap K_{H}\right] \subseteq
\partial _{m}^{v}\left( NG_{n,m} \cap D_{n,m}\right)
\end{equation*}%
where
$$
K_{H}=\bigcap\limits_{i=0}^{n-1}\ker d_i^{h^{(n-1,m)}}
$$
and
\begin{equation*}
K_{I}=\bigcap\limits_{i\in I}\ker d_{i}^{^{v}(n,m-1)}\ \ \
\text{and}\ \ K_{J}=\bigcap\limits_{j\in J}\ker d_{j}^{^{v}(n,m-1)}.
\end{equation*}%
Similarly, for $n\geqslant 2, m\geqslant 1$ and $I', J'\subseteq [n-1] $ with
$I'\cup J'=[n-1]$, there is the inclusion
\begin{equation*}
\left[ K_{I'}\cap K_{V},K_{J'}\cap K_{V}\right] \subseteq
\partial _{n}^{h}\left( NG_{n,m}\cap D_{n,m}\right)
\end{equation*}%
where
$$
K_{V}=\bigcap\limits_{i=0}^{m-1}\ker d_i^{v^{(n,m-1)}}
$$
and
\begin{equation*}
K_{I'}=\bigcap\limits_{i\in I'}\ker d_{i}^{^{h}(n-1,m)}\ \ \
\text{and}\ \ K_{J'}=\bigcap\limits_{j\in J'}\ker
d_{j}^{^{h}(n-1,m)}.
\end{equation*}
\end{prop}

\begin{proof}
 We know that from the results of \cite{mutpor2}, there are already the following inclusions in both horizontal and vertical directions: if $m$ is constant in the horizontal direction, then for $n\geqslant 2$, we have
$$ \left[ K_{I'},K_{J'}\right] \subseteq
\partial _{n}^{h}\left( NG_{n,m}\cap D_{n,m}\right)
 $$
 and if $n$ is constant in the vertical  direction, then  for $m \geqslant 2$, we have
 $$
 \left[ K_{I},K_{J}\right] \subseteq
\partial _{m}^{v}\left( NG_{n,m} \cap D_{n,m}\right).
 $$
Thus the  result can be seen easily by using these inclusions.
 \end{proof}

Now, in low dimensions, we investigate the images of these functions under the boundary homomorphisms $\partial_n^h$ and $\partial_m^v$.

For  $NG_{1,2}$ take $x,y\in NG_{1,1}=\ker {d_{0}^{{h}^{(1,1)}}}\cap \ker ${$%
d_{0}^{v^{(1,1)}}$}, then we obtain
\begin{eqnarray*}
\partial_{2}^{^{v}(12)}(F_{(\left( \emptyset ,(0)\right)
,\left( \emptyset ,(1)\right) })(x,y) &=&d_{2}^{v^{(12)}}\left( [ s_{0}^{v^{(11)}}\left( x\right)
,s_{1}^{v^{(11)}}(y)][ s_{1}^{v^{(11)}}y ,s_{1}^{v^{(11)}}x] \right) \\
&=&[s_{0}^{v^{(01)}}d_{1}^{v^{(11)}}x, y ][y,x]
\end{eqnarray*}
where $[s_{0}^{v^{(01)}}d_{1}^{v^{(11)}}x, y ][y,x]\in \lbrack \ker d_{0}^{^{v}(11)},\ker d_{1}^{^{v}(11)}]$ from \cite{mutpor2}. Further we obtain
$$
d_{0}^{h^{(11)}}([s_{0}^{v^{(01)}}d_{1}^{v^{(11)}}x, y ][y,x])=1,$$
hence $[s_{0}^{v^{(01)}}d_{1}^{v^{(11)}}x, y ][y,x]\in [\ker d_{0}^{{v}^{(11)}}\cap \ker
d_{0}^{h^{(11)}},\ker d_{1}^{v^{(11)}}\cap \ker
d_{0}^{h^{(11)}}].$

Similarly  for $y\in NG_{0,2}$ and $x\in NG_{1,1}$, from $NG_{1,1}\times NG_{0,2}$ to $NG_{1,2}$
we obtain
\begin{eqnarray*}
\partial_{2}^{v^{(12)}}(F_{(\left( \emptyset ,(1)\right)
,\left( (0) ,\emptyset \right) })(x,y) &=&d_{2}^{v^{(12)}}[s_1^{v^{(11)}}(x),s_0^{h^{(02)}}(y)]\\
&=&[x,d_{2}^{v^{(12)}}s_0^{h^{(02)}}(y)].
\end{eqnarray*}

Since $d_0^{h^{(11)}}[x,d_{2}^{v^{(12)}}s_0^h(y)]=1$ we have $[x,d_{2}^{v^{(12)}}s_0^h(y)] \in \ker d_0^{h^{(11)}}.$ Furthermore,
$[x,d_{2}^{v^{(12)}}s_0^h(y)] \in [\ker d_0^{v^{(11)}} ,\ker d_1^{v^{(11)}}]$.

By a similar way one can  show that the images of other generating elements in
$$[\ker d_{0}^{v^{(11)}}\cap \ker d_{0}^{h^{(11)}},\ker d_{1}^{v^{(11)}}\cap \ker d_{0}^{h^{(11)}}].$$

Thus we have the following equality%
\begin{equation*}
\partial _{2}^{^{v}}\left( NG_{1,2}\cap D_{1,2}\right) =[\ker d_{0}^{^{v}(11)}\cap \ker
d_{0}^{^{h}(11)},\ker
d_{1}^{^{v}(11)}\cap \ker
d_{0}^{^{h}(11)}].
\end{equation*}

For $NG_{2,1}$ take $x,y\in NG_{1,1}$. Then we obtain
\begin{eqnarray*}
\partial_{2}^{^{h}(2,1)}(F_{(((0),\emptyset ),\left(
(1),\emptyset \right) )})(x,y)=[s_{0}^{h^{(01)}}d_{1}^{h^{(11)}}(x),y][y,x]
\end{eqnarray*}
where $[s_{0}^{h^{(01)}}d_{1}^{h^{(11)}}(x),y][y,x]\in [
\ker d_{0}^{^{h}(11)},\ker d_{1}^{^{h}(11)}]$ and $$%
d_{0}^{^{v}(11)}([s_{0}^{h^{(0,1)}}d_{1}^{h^{(1,1)}}(x),y][y,x])=1,$$
hence
$$[s_{0}^{h^{(0,1)}}d_{1}^{h^{(1,1)}}(x),y][y,x]\in \ker
d_{0}^{^{v}(11)}\cap [ \ker d_{0}^{^{h}(11)},\ker d_{1}^{^{h}(11)}].$$

Thus we obtain the following equality%
\begin{equation*}
\partial _{2}^{h}\left( NG_{2,1}\cap D_{2,1}\right) =[ \ker d_{0}^{^{h}(11)}\cap \ker
d_{0}^{^{v}(11)},\ker
d_{1}^{^{h}(11)}\cap \ker
d_{0}^{^{v}(11)}].
\end{equation*}
For \ $(n,m)=(2,1)$ and $(n,m)=(1,2)$, we can summarize these situations in the following diagram.
$$
\begin{tabular}{|c|c|c|c|c|}
\hline
$\underline{\alpha}$& $\underline{\beta}$&$I'$&$J'$& $V$\\
\hline
$\left( (0), \emptyset \right) $ & $\left((1), \emptyset \right) $ & $ \{0\}$&$\{1\}$&$\{0\}$ \\
\hline
$\underline{\alpha}$& $\underline{\beta}$&$I$&$J$& $H$\\
\hline
$\left( \emptyset,(0) \right) $ & $\left( \emptyset ,(1)\right) $ & $ \{0\}$&$\{1\}$&$\{0\}$ \\
\hline
\end{tabular}
$$

For $NG_{2,2}$ take $x,y\in NG_{2,1}=\ker {d_{0}^{{h}^{(2,1)}}}\cap \ker {%
d_{1}^{{h}^{(2,1)}}}\cap \ker ${$d_{0}^{v^{(2,1)}}$}, then we obtain
\begin{eqnarray*}
\partial _{2}^{^{v}(22)}(F_{(\left( \emptyset ,(0)\right) ,\left( \emptyset
,(1)\right) })(x,y) &=&d_{2}^{v^{(22)}}\left(
[s_{0}^{v^{(21)}}(x),s_{1}^{v^{(21)}}(y)][s_{1}^{v^{(21)}}(y),s_{1}^{v^{(21)}}(x)]]\right)
\\
&=&[s_{0}^{v^{(20)}}d_{1}^{v^{(21)}}x,y][y,x]
\end{eqnarray*}%
where $[s_{0}^{v^{(20)}}d_{1}^{v^{(21)}}x,y][y,x]\in \lbrack \ker
d_{0}^{^{v}(21)},\ker d_{1}^{^{v}(21)}]$ from \cite{mutpor2}. Further we
obtain
\begin{equation*}
d_{0}^{h^{(21)}}([s_{0}^{v^{(20)}}d_{1}^{v^{(21)}}x,y][y,x])=1,
\end{equation*}%
hence $[s_{0}^{v^{(20)}}d_{1}^{v^{(21)}}x,y][y,x]\in \lbrack \ker d_{0}^{{v}%
^{(21)}}\cap \ker d_{0}^{h^{(21)}},\ker d_{1}^{v^{(21)}}\cap \ker
d_{0}^{h^{(21)}}].$

Similarly for  $x\in NG_{1,2}\ $and $y\in NG_{2,1}$, from $NG_{1,2}\times
NG_{2,1}$ to $NG_{2,2}$ we obtain
\begin{eqnarray*}
\partial _{2}^{v^{(22)}}(F_{((1),\emptyset ),(\emptyset ,(0))}(x,y))
&=&d_{2}^{v^{(22)}}[s_{1}^{h^{(12)}}(x),s_{0}^{v^{(21)}}(y)] \\
&=&[d_{2}^{v^{(22)}}s_{1}^{h^{(12)}}(x),d_{2}^{v^{(22)}}s_{0}^{v^{(21)}}(y)].
\end{eqnarray*}

Since $d_{0}^{h^{(11)}}[x,d_{2}^{v^{(12)}}s_{0}^{h}(y)]=1$ we have $%
[x,d_{2}^{v^{(12)}}s_{0}^{h}(y)]\in \ker d_{0}^{h^{(11)}}.$ Furthermore, $%
[x,d_{2}^{v^{(12)}}s_{0}^{h}(y)]\in \lbrack \ker d_{0}^{v^{(11)}},\ker
d_{1}^{v^{(11)}}]$.

By a similar way, one can show that the images of other generating elements
in
\begin{equation*}
\lbrack \ker d_{0}^{v^{(21)}}\cap \ker d_{0}^{h^{(21)}},\ker
d_{1}^{v^{(21)}}\cap \ker d_{0}^{h^{(21)}}].
\end{equation*}

Thus we have the following equality%
\begin{equation*}
\partial _{2}^{^{v}}\left( NG_{2,2}\cap D_{2,2}\right) =[\ker
d_{0}^{^{v}(21)}\cap \ker d_{0}^{^{h}(21)},\ker d_{1}^{^{v}(21)}\cap \ker
d_{0}^{^{h}(21)}].
\end{equation*}%
We can summarize this in the following diagram for $(n,m)=(2,2)$.
$$
\begin{tabular}{|c|c|c|c|c|}
\hline
$\underline{\alpha}$& $\underline{\beta}$&$I'$&$J'$& $V$\\
\hline
$\left( (0), \emptyset \right) $ & $\left((1), \emptyset \right) $ & $ \{0\}$&$\{1\}$&$\{0,1\}$ \\
\hline
$\underline{\alpha}$& $\underline{\beta}$&$I$&$J$& $H$\\
\hline
$\left( \emptyset,(0) \right) $ & $\left( \emptyset ,(1)\right) $ & $ \{0\}$&$\{1\}$&$\{0,1\}$ \\
\hline
\end{tabular}
$$
Using the calculation method given above we obtained the following equalities in low dimensions.
\begin{equation*}
\begin{tabular}{lll}
$\partial_{2}^{^{v}}\left( NG_{0,2}\cap D_{0,2}\right) =[\ker d_{0}^{v^{(01)}}, \ker d_{1}^{v^{(01)}}],$  \\
$\partial _{2}^{^{v}}\left( NG_{1,2}\cap
D_{1,2}\right) =[\ker d_{0}^{^{v}(11)}\cap \ker d_{0}^{^{h}(11)},\ker
d_{1}^{^{v}(11)}\cap \ker d_{0}^{^{h}(11)}] $, \\
$ \partial
_{2}^{^{v}}\left( NG_{2,2}\cap D_{2,2}\right) = \lbrack \ker
d_{0}^{^{v}(21)}\cap \ker d_{0}^{^{h}(21)}\cap \ker d_{1}^{^{h}(21)},\ker d_{1}^{^{v}(21)}\cap \ker d_{0}^{^{h}(21)}\cap \ker d_{1}^{^{h}(21)}]$,  \\
$ \partial
_{2}^{^{h}}\left( NG_{2,0}\cap D_{2,0}\right)=[\ker d_{0}^{^{h}(10)}, \ker d_{1}^{^{h}(10)}] $,   \\
$ \partial _{2}^{^{h}}\left( NG_{2,1}\cap
D_{2,1}\right)=\lbrack \ker d_{0}^{^{h}(11)}\cap \ker d_{0}^{^{v}(11)},\ker
d_{1}^{^{h}(11)}\cap \ker d_{0}^{^{v}(11)}\rbrack $ ,  \\
$\partial
_{2}^{^{h}}\left( NG_{2,2}\cap D_{2,2}\right)=\lbrack \ker
d_{0}^{^{h}(12)}\cap \ker d_{0}^{^{v}(12)}\cap \ker d_{1}^{^{v}(12)},\ker d_{1}^{^{h}(12)}\cap \ker d_{0}^{^{v}(12)}\cap \ker d_{1}^{^{v}(12)}\rbrack$ .
\end{tabular}%
\end{equation*}

\section{Applications}

\subsection{ Crossed modules from bisimplicial groups}

Recall that a
 \textit{crossed module} is a group morphism $\partial:M \rightarrow
P$ endowed with a (left) action of $P$ on $M$ such that

$CM1.$ the morphism $\partial$ is $P$-equivariant, where $P$ acts on
itself via conjugacy,

$CM2.$ for $x,y\in M$ we have $^{\partial x}y=xyx^{-1}$.

\begin{prop}\label{2}
Let $\mathbf{G}_{\ast, \ast}$ be a bisimplicial group with Moore
bicomplex $\mathbf{NG}_{\ast, \ast}$. If $p\geqslant  1$ and
$q\geqslant  1$, $NG_{p,q}=\{1\}$, then the maps $d_{1}^v{^{(01)}}$, $d_{1}^h{^{(10)}}$ and
\begin{eqnarray*}
\partial :NG_{0,1}\times NG_{1,0}&\rightarrow &NG_{0,0}\\
(x,y)&\mapsto & d_{1}^v{^{(01)}}(x)d_{1}^h{^{(01)}}(y)
\end{eqnarray*}
are crossed modules.
\end{prop}
\begin{proof}
See Appendix B
\end{proof}

\subsection{Crossed squares from bisimplicial groups}

First, we recall  detailed definition of a crossed square from \cite{loday}.

A \emph{crossed square} of groups is a commutative square of
group morphisms
\begin{equation*}
\xymatrix{ L \ar[d]_{\lambda'} \ar[r]^{\lambda} & M \ar[d]^{\mu}
\\ N \ar[r]_{\nu} & P}
\end{equation*}%
with action of $P$ on every other group and a map $h:M\times N \rightarrow L$ such that

\begin{enumerate}
\item  The maps $\lambda$ and $\lambda'$ are $P$-equivariant and $\nu$,$\mu$, $\mu \circ \lambda$ and $\nu \circ \lambda'$ are crossed modules,

\item $\lambda\circ h(x,y)=x ^{\nu(y)}x^{-1}, \ \lambda'\circ h(x,y)=(^{\mu(x)}y)y^{-1},$

\item  $h(\lambda(z),y)=z^{\nu(y)}z^{-1}, \ h(x,\lambda'(z))=(^{\mu(x)}z)z^{-1}$,

\item $h(xx',y)=^{\mu(x)}h(x',y)h(x,y)$, \ $h(x,yy')=h(x,y)^{\nu(y)}h(x,y')$,

\item $h(^tx,^ty)=^th(x,y)$
\end{enumerate}
for $x,x'\in M$, $y,y'\in N, z\in L$ and $t\in P$.

The following proposition was initially given by Conduch\'e in \cite{cond}. Here, we see that the $h$-map of the crossed square can be given by the function
$F_{(\emptyset,(0)),((0),\emptyset)}:NG_{01}\times NG_{10}\rightarrow NG_{11}$.

\begin{prop}\label{3}
 Let $\mathbf{G}_{\ast,\ast}$ be a bisimplicial group and $\mathbf{NG}_{\ast,\ast}$ its Moore
bicomplex. Suppose $NG_{k_1,k_2}=\{1\}$  for any $k_1\geqslant  2$ or
$k_2\geqslant 2$. Then the diagram
$$\xymatrix{
 & NG_{1,1} \ar[d]_{\partial_1^{{v}^{(11)}}} \ar[r]^{\partial_1^{{h}^{(11)}}}
              &   NG_{0,1}\ar[d]^{\partial_1^{{v}^{(01)}}}
 \\ & NG_{1,0}  \ar[r]_{\partial_1^{{h}^{(10)}}}
               &  NG_{0,0}      }$$
is a crossed square. $NG_{0,0}$ acts on other groups via the
degeneracies $s_{0}^{h}$ and $s_{0}^{v}. $ The $h$-map is given by
the map $ F_{(\emptyset,(0)),((0),\emptyset)}(x,y) $, namely,
\begin{eqnarray*}
h:NG_{0,1}\times NG_{1,0}&\rightarrow& NG_{1,1}\\
 (x,y)&\mapsto &h(x,y)=F_{(\emptyset,(0)),((0),\emptyset)}(x,y)
\end{eqnarray*}
for $ x\in NG_{0,1},\ y\in NG_{1,0}$ where $(\emptyset,(0)),((0),\emptyset)\in S(1)\times S(1)$ and
$$
F_{(\emptyset,(0)),((0),\emptyset)}(x,y)=[s_{0}^{h^{(01)}}(x),s_{0}^{v^{(10)}}(y)].
$$
\end{prop}
\begin{proof}
For the axioms  see Appendix B.
\end{proof}

\textbf{Remark:} This result can be extended to crossed $n$-cubes defined by Ellis and Steiner in \cite{els}.
In this case, if  $\mathbf{G}_{\ast _{1}\ast _{2}\dots\ast _{n}}$ is an $n$-simplicial group with Moore $n$-complex $\mathbf{NG}_{\ast _{1}\ast
_{2}\dots\ast _{n}},$ such that $\ NG_{\ast _{1}\ast _{2}\dots\ast
_{n}}=\left\{ 1\right\}$ for any $\ast _{j}\geqslant  2,\ (1\leqslant j\leqslant n),$ then we can easily say that, as a generalization of the above result, this Moore $n$-complex has a crossed $n$-cube structure. Then the $h$-maps of associated crossed $n$-cube are given by the functions $F_{\alpha,\beta}$ in the Moore $n$-complex.

\subsection{2-Crossed modules from bisimplicial groups}

To characterize the simplicial groups having a Moore complex trivial in dimension bigger or equal to three, Conduch\'e has \cite{con}. The definition of
a 2-crossed module is also recalled for completeness \cite{con}.

\emph{A 2-crossed module} is a complex of length 2
$$
\xymatrix{L\ar[r]^-{\partial_2}&M\ar[r]^{\partial_1}&N}
$$
of $N$-groups with action of $M$ on $L$ and a function $\{,\}:M\times M \rightarrow L$, called `Peiffer lifting', such that
\begin{enumerate}
\item $ \partial_2\{y,y'\}=yy'y^{-1}(^{\partial_1(y)}y')^{-1}$,
\item $\{\partial_2z,\partial_2 z'\}=zz'z^{-1}z'^{-1},$
\item $\{\partial_2z,y\}\{y,\partial_2z\}=z(^{\partial_1y}z)^{-1}, $
\item $\{y,y'y''\}=\{y,y'\}\{y,y''\}\{\partial_2\{y,y''\}^{-1},^{\partial_1y}y'\},$\
\item $\{yy',y''\}=\{y,y'y''y'^{-1}\}^{\partial_1y}\{y',y''\},$
\item $^{x}\{y,y'\}=\{^xy,^xy'\}$
\end{enumerate}
for $x\in N,y,y',y''\in M$ and $z,z'\in L$.

Conduch\'e  proved the following proposition in \cite{con}.

\begin{proposition} \label{4}
Let $\mathbf{G}_{\ast}$ be a simplicial group and $\mathbf{NG}_*$ its Moore complex. Suppose $NG_n=\{1\}$ for $n\geqslant 3.$ Then
the complex
$$
\xymatrix{NG_2\ar[r]^-{\partial_2}&NG_1\ar[r]^{\partial_1}&NG_0}
$$
is a 2-crossed module, where $NG_0$ acts on $NG_1$ and $NG_2$ by conjugacy via the degeneracies and the Peiffer lifting is given by
$$
\{y,y'\}=s_1(yy'y^{-1})s_0(y)s_1(y')^{-1}s_0(y)^{-1}.
$$
\end{proposition}
Conduch{\'{e}} also constructed in \cite{cond} a $2$-crossed module from a crossed square
\begin{equation*}
\xymatrix{ L \ar[d]_{\lambda'} \ar[r]^{\lambda} & M \ar[d]^{\mu} \\ N
\ar[r]_{\nu} & P }
\end{equation*}%
as
\begin{equation*}
L\llabto{2}{(\lambda ^{-1},\lambda ^{\prime })}M\rtimes N
\llabto{2}{\mu \nu}P.
\end{equation*}
Since
$$\xymatrix{
  NG_{1,1} \ar[d]_{\partial_1^{{v}}} \ar[r]^{\partial_1^{{h}}}
               & NG_{0,1}\ar[d]^{\partial_1^{{v}}}\\
 NG_{1,0}  \ar[r]_{\partial_1^{{h}}}
               &  NG_{0,0}      }$$
is a crossed square, the complex of morphisms of groups
\begin{equation*}
NG_{1,1}\llabto{2}{((d_1^h)^{-1}, d_1^{v})} NG_{0,1}\times NG_{1,0}
\llabto{2}{d_1^{v^{(01)}}d_1^{h^{(10)}}}NG_{0,0}.
\end{equation*}
is a 2-crossed module with the Peiffer lifting map
$$
\{(x,a),(y,b)\}=F_{(\emptyset,(0)),((0),\emptyset)}(x,ab)
$$
for $(x,a),(y,b)\in NG_{0,1}\times NG_{1,0}$ as given in \cite{cond}.

We get the following result.
\begin{proposition}\label{5}
Let $\mathbf{G}_{\ast,\ast}$ be a bisimplicial group with Moore bicomplex $\mathbf{NG}_{\ast,\ast}$. If for any $p\geqslant 0$ and $q\geqslant 3,$
$NG_{p,q}=\{1\},$  then the complex of morphisms of groups
$$
\xymatrix{NG_{p,2}\ar[r]^-{\partial_2^v}&NG_{p,1}\ar[r]^{\partial_1^v}&NG_{p,0}}
$$
is a 2-crossed module. The Peiffer lifting map
$$
\{-,-\}:NG_{p,1}\times NG_{p,1}\to NG_{p,2}
$$
is given  by
$$ \{x,y\}=\left(F_{(\emptyset,(0)),(\emptyset,(1))}(x,y)\right)^{-1}=[s_{1}^{v^{(p1)}}(x),s_{1}^{v^{(p1)}}(y)][s_{1}^{v^{(p1)}}(y),s_{0}^{v^{(p1)}}(x)] $$
for $ x,y\in NG_{p,1} $. Similarly, if $NG_{p,q}=\{1\}$ for any $q\geqslant 0$ and $p\geqslant 3,$
then the complex of morphisms of groups
$$
\xymatrix{NG_{2,q}\ar[r]^-{\partial_2^h}&NG_{1,q}\ar[r]^{\partial_1^h}&NG_{0,q}}
$$
is a 2-crossed module. The Peiffer lifting map
$$
\{-,-\}:NG_{1,q}\times NG_{1,q}\to NG_{2,q}
$$
is given by
$$ \{x,y\}=\left(F_{((0),\emptyset),((1),\emptyset)}(x,y)\right)^{-1}=[s_{1}^{h^{(1q)}}(x),s_{1}^{h^{(1q)}}(y)][s_{1}^{h^{(1q)}}(y),s_{0}^{h^{(1q)}}(x)]$$
for  $ x,y\in NG_{1,q} $.
\end{proposition}
\begin{proof}
See Appendix B.
\end{proof}

\section{Appendices}
\subsection{Appendix A}
For $(n,m)=(1,1)$, consider the set
$$S(1)\times S(1)=\{(\emptyset,\emptyset),(\emptyset,(0)),((0),\emptyset),((0),(0))\}.$$

\begin{enumerate}
\item Take $\underline{\alpha}=(\emptyset,\emptyset)$ and $\underline{\beta}=(\emptyset,(0))$. In this case, the function $F_{\underline{\alpha},\underline{\beta}}$ becomes from $NG_{1,1}\times NG_{1,0}$ to $NG_{1,1}.$ This map can be defined for any $x\in NG_{1,1}$ and $y\in NG_{1,0}$ by
$$
F_{(\emptyset,\emptyset),(\emptyset,(0))}(x,y)=[x,s_0^{v^{(10)}}(y)].
$$

\item Take $\underline{\alpha}=(\emptyset,\emptyset)$ and $\underline{\beta}=((0), \emptyset)$. In this case, the function $F_{\underline{\alpha},\underline{\beta}}$ becomes
from $NG_{1,1}\times NG_{0,1}$ to $NG_{1,1}.$ This map can be defined by
$$
F_{(\emptyset,\emptyset),((0),\emptyset)}(x, y)=[x,s_0^{h^{(01)}}(y)].
$$
for any $x\in NG_{1,1}$ and $y\in NG_{0,1}$.

\item For $\underline{\alpha}=(\emptyset,\emptyset)$ and $\underline{\beta}=((0),(0))$. The map
$$
F_{(\emptyset,\emptyset),((0),(0))}:NG_{1,1}\times NG_{0,0}\rightarrow NG_{1,1}
$$
is defined by
$$
F_{(\emptyset,\emptyset),((0),(0))}(x, y)=[x,(s_0^{v^{(00)}}s_0^{h^{(01)}}(y))]
$$
for all $x\in NG_{1,1}$ and $y\in NG_{00}$.

\item For $\underline{\alpha}=((0),\emptyset)$ and $\underline{\beta}=(\emptyset,(0))$. Then the map
$$
F_{((0),\emptyset),(\emptyset,(0))}:NG_{0,1}\times NG_{1,0}\rightarrow NG_{1,1}
$$
can be calculated for any $x\in NG_{0,1}$ and $y\in NG_{1,0}$ by
\begin{align*}
F_{\left( (0),\emptyset \right) ,\left( \emptyset ,(0)\right) }(x,y)  =&p\mu
(s_{0}^{h^{(01)}}(x),s_{0}^{v^{(10)}}(y)) \\
=&p_{0}^{h}p_{0}^{v}(s_{0}^{h^{(01)}}(x)s_{0}^{v^{(10)}}(y)s_{0}^{h^{(01)}}(x)^{-1}s_{0}^{v^{(10)}}(y)^{-1}) \\
=&p_{0}^{h}((s_{0}^{v^{(10)}}(y)s_{0}^{h^{(01)}}(x)s_{0}^{v^{(10)}}(y)^{-1}s_{0}^{h^{(01)}}(x)^{-1})\\
&s_{0}^{v^{(10)}}d_{0}^{v^{(11)}}(s_{0}^{h^{(01)}}(x)s_{0}^{v^{(10)}}(y)s_{0}^{h^{(01)}}(x)^{-1}s_{0}^{v^{(10)}}(y)^{-1})) \\
=&p_{0}^{h}(s_{0}^{v^{(10)}}(y)s_{0}^{h^{(01)}}(x)s_{0}^{v^{(10)}}(y)^{-1}s_{0}^{h^{(01)}}(x)^{-1})(s_{0}^{v^{(10)}}d_{0}^{v^{(11)}}
\\
&s_{0}^{h^{(01)}}(x)s_{0}^{v^{(10)}}d_{0}^{v^{(11)}}s_{0}^{v^{(10)}}(y)s_{0}^{v^{(10)}}d_{0}^{v^{(11)}}s_{0}^{h^{(01)}}(x)^{-1}s_{0}^{v^{(10)}}d_{0}^{v^{(11)}}s_{0}^{v^{(10)}}(y)^{-1}) \\
=&p_{0}^{h}(s_{0}^{v^{(10)}}(y)s_{0}^{h^{(01)}}(x)s_{0}^{v^{(10)}}(y)^{-1}s_{0}^{h^{(01)}}(x)^{-1})
(s_{0}^{v^{(10)}}s_{0}^{h^{(00)}}d_{0}^{v^{(01)}}(x)\\
&s_{0}^{v^{(10)}}(y)s_{0}^{v^{(10)}}s_{0}^{h^{(00)}}d_{0}^{v^{(01)}}(x)^{-1}s_{0}^{v^{(10)}}(y)^{-1}) \\
 =&p_{0}^{h}(s_{0}^{v^{(10)}}(y)s_{0}^{h^{(01)}}(x)s_{0}^{v^{(10)}}(y)^{-1}s_{0}^{h^{(01)}}(x)^{-1}) \\
 =&(s_{0}^{h^{(01)}}(x)s_{0}^{v^{(10)}}(y)s_{0}^{h^{(01)}}(x)^{-1}s_{0}^{v^{(10)}}(y)^{-1})
\\
&s_{0}^{h^{(01)}}d_{0}^{h^{(11)}}(s_{0}^{v^{(10)}}(y)s_{0}^{h^{(01)}}(x)s_{0}^{v^{(10)}}(y)^{-1}s_{0}^{h^{(01)}}(x)^{-1}) \\
 =&(s_{0}^{h^{(01)}}(x)s_{0}^{v^{(10)}}(y)s_{0}^{h^{(01)}}(x)^{-1}s_{0}^{v^{(10)}}(y)^{-1}) \\
  =&[s_{0}^{h^{(01)}}(x),s_{0}^{v^{(10)}}(y)].
\end{align*}

\item For $\underline{\alpha}=((0),\emptyset)$ and $\underline{\beta}=((0),(0))$. Then the map
$$
F_{((0),\emptyset),((0),(0))}:NG_{0,1}\times NG_{0,0}\rightarrow NG_{1,1}
$$
can be calculated for any $x\in NG_{0,1}$ and $y\in NG_{0,0}$ by
\begin{align*}
F_{((0),\emptyset),((0),(0))}(x, y)&=p\mu(s_{\underline{\alpha}}, s_{\underline{\beta}})(x, y)\\
&=p_0^hp_0^v[s_0^{h^{(01)}}(x),s_0^{h^{(01)}}s_0^{v^{(00)}}(y)]\\
&=1
\end{align*}

\item Similarly for $\underline{\alpha}=(\emptyset,(0))$ and $\underline{\beta}=((0),(0))$, the map
$$
F_{(\emptyset, (0)),((0),(0))}:NG_{1,0}\times NG_{0,0}\rightarrow NG_{1,1}
$$
is the identity as given in the previous step.
\end{enumerate}

 By taking $(n,m)$= $(0,2)$ and $(2,0)$, we calculate the possible non identity maps with codomain $NG_{0,2}$ and $ NG_{2,0}$ respectively.

First $(n,m)=(0,2)$. Consider the set
$$
S(0)\times S(2)=\{(\emptyset,\emptyset),(\emptyset,(0)),(\emptyset,(1)),(\emptyset,(1,0))\}.
$$
We try to find the functions $F_{\underline{\alpha},\underline{\beta}}$ with codomain $NG_{0,2}$. In this case the only non identity map $F_{\underline{\alpha},\underline{\beta}}$ can be defined by choosing $\underline{\alpha}=(\emptyset,(0))$ and $\underline{\beta}=(\emptyset,(1))$. Then this  is a map from $NG_{0,1}\times NG_{0,1}$ to $NG_{0,2}$. This map is calculated as follows. For $x,y\in NG_{0,1}$,
we obtain
\begin{align*}
F_{(\emptyset,(0)),(\emptyset,(1))}(x, y)&=p\mu(s_{\underline{\alpha}}, s_{\underline{\beta}})(x, y)\\
&=p_1^vp_0^v[s_0^{v^{(01)}}(x),s_1^{v^{(01)}}(y)]\\
&=[ s_{0}^{v^{(01)}}x ,s_{1}^{v^{(01)}}y] [s_{1}^{v^{(01)}}y ,s_{1}^{h^{(01)}}x ]
\in NG_{0,2}.
\end{align*}

Now suppose $(n,m)=(2,0)$. From the set
$$
S(2)\times S(0)=\{(\emptyset,\emptyset),((0),\emptyset),((1),\emptyset),((1,0),\emptyset)\}
$$
we can choose $\underline{\alpha}=((0),\emptyset)$ and $\underline{\beta}=((1),\emptyset)$. Then  $F_{\underline{\alpha},\underline{\beta}}$ is a map from $NG_{1,0}\times NG_{1,0}$ to $NG_{2,0}$. This map can be given by for $x,y\in NG_{1,0}$
\begin{align*}
F_{((0),\emptyset),((1),\emptyset)}(x, y)&=p\mu(s_{\underline{\alpha}},s_{\underline{\beta}})(x,y)\\
&=p_1^hp_0^h[s_0^{h^{(10)}}(x),s_1^{h^{(10)}}(y)]\\
&=[s_0^{h^{(10)}}(x),s_1^{h^{(10)}}(y)][s_1^{h^{(10)}}(y),s_1^{h^{(10)}}(x)]\in NG_{2,0}.
\end{align*}

Now, by taking $(n,m)$= $(1,2)$ and $(2,1)$, we shall define the possible non identity maps $F_{\underline{\alpha},\underline{\beta}}$ whose codomain $NG_{1,2}$ and $ NG_{2,1}$ respectively.

First suppose that $(n,m)=(1,2)$. We set
$$
S(1)\times S(2)=\{(\emptyset,\emptyset),(\emptyset,(1)),(\emptyset,(0)),(\emptyset,(1,0)),((0),\emptyset),((0),(1)),((0),(0)),((0),(1,0))\}.
$$
In the following calculations, by taking appropriate $\underline{\alpha},\underline{\beta}$ from the set $S(1)\times S(2)$, we shall give all the non identity maps whose codomain $NG_{1,2}$. To obtain these maps, we can choose the possible $\underline{\alpha},\underline{\beta}$ from the set $S(1)\times S(2)$ as follows.

\begin{center}
\begin{tabular}{ll}
1. $(\underline{\alpha},\underline{\beta})=((\emptyset,(0)),(\emptyset,(1)))$ &
2. $(\underline{\alpha},\underline{\beta})=((\emptyset,(1)),((0),\emptyset))$\\
3. $(\underline{\alpha},\underline{\beta})=((\emptyset,(0)),((0),\emptyset))$ & 4. $(\underline{\alpha},\underline{\beta})=(((0),(1)),(\emptyset,(0)))$\\
5. $(\underline{\alpha},\underline{\beta})=(((0),(0)),(\emptyset,(1))).$ &
\end{tabular}
\end{center}
Now we calculate the functions $F_{\underline{\alpha},\underline{\beta}}$ for these pairings $(\underline{\alpha},\underline{\beta})$.
\begin{enumerate}
\item For $\underline{\alpha}=(\emptyset,(0))$ and $\underline{\beta}=(\emptyset,(1))$, we obtain the map
 $$
 F_{(\emptyset,(0)),(\emptyset,(1))}:NG_{1,1}\times NG_{1,1} \longrightarrow NG_{1,2}.
 $$
This map can be given by
\begin{align*}
F_{(\emptyset, (0)),(\emptyset, (1))}(x, y)=&p\mu(s_{\underline{\alpha}},s_{\underline{\beta}})(x,y)\\
=& p_1^v p_0^v p_0^h[s_0^{v^{(11)}}x,s_1^{v^{(11)}}y]\\
=&[s_0^{v^{(11)}}(x),s_1^{v^{(11)}}(y)][s_1^{v^{(11)}}(y),s_1^{v^{(11)}}(x)]\in NG_{1,2}
\end{align*}
for $x,y\in NG_{1,1}$.

\item For $\underline{\alpha}=(\emptyset,(1)), \underline{\beta}=((0),\emptyset)$, we get the map
$$
 F_{(\emptyset,(1)),((0),\emptyset)}:NG_{1,1}\times NG_{0,2} \longrightarrow NG_{1,2}
$$
given by
$$
F_{(\emptyset, (1)),((0),\emptyset)}(x, a)=[s_1^{v^{(11)}}(x),s_0^{h^{(02)}}(a)]\in NG_{1,2}
$$
for $x\in NG_{1,1}$ and $a\in NG_{0,2}$.

\item For $\underline{\alpha}=(\emptyset,(0)), \underline{\beta}=((0),\emptyset)$, we have the following map
$$
 F_{(\emptyset,(0)),((0),\emptyset)}:NG_{1,1}\times NG_{0,2} \longrightarrow NG_{1,2}
$$
calculated by
$$
F_{(\emptyset, (0)),((0),\emptyset)}(x, a)=[s_0^{v^{(11)}}(x),s_0^{h^{(02)}}(a)]\in NG_{1,2}
$$
for $x\in NG_{1,1}$ and $a\in NG_{0,2}$.

\item For $\underline{\alpha}=((0),(1)), \underline{\beta}=(\emptyset,(0))$, we get the following map
$$
F_{(0),(1)),(\emptyset,(0))}:NG_{0,1}\times NG_{1,1} \longrightarrow NG_{1,2}
$$
given by
$$
F_{((0), (1)),(\emptyset,(0))}(x, y)=[s_0^{h^{(02)}}s_1^{v^{(01)}}(x),s_0^{v^{(11)}}(y)]\in NG_{1,2}
$$
for $x\in NG_{0,1}$ and $y\in NG_{1,1}$.

\item For $\underline{\alpha}=((0),(0)), \underline{\beta}=(\emptyset,(1))$, we get the following map
$$
F_{(0),(0)),(\emptyset,(1))}:NG_{0,1}\times NG_{1,1} \longrightarrow NG_{1,2}
$$
given by
$$
F_{((0), (0)),(\emptyset,(1))}(x, y)=[s_0^{h^{(02)}}s_0^{v^{(01)}}(x),s_1^{v^{(11)}}(y)]\in NG_{1,2}
$$
for $x\in NG_{0,1}$ and $y\in NG_{1,1}$.
\end{enumerate}

Now suppose that $(n,m)=(2,1)$. We consider the set $S(2)\times S(1)$. By choosing appropriate $\underline{\alpha},\underline{\beta}$ from the set $S(2)\times S(1)$, we can calculate similarly all the non identity maps with codomain $NG_{2,1}$. To obtain these maps, we take the possible $\underline{\alpha},\underline{\beta}$ as follows.
\begin{center}
\begin{tabular}{ll}
1. $(\underline{\alpha},\underline{\beta})=(((0),\emptyset),((1),\emptyset,))$ &
2. $(\underline{\alpha},\underline{\beta})=(((1),\emptyset),(\emptyset,(0)))$\\
3. $(\underline{\alpha},\underline{\beta})=(((0),\emptyset),(\emptyset,(0)))$ &
4. $(\underline{\alpha},\underline{\beta})=(((1),(0)),((0),\emptyset))$\\
5. $(\underline{\alpha},\underline{\beta})=(((0),(0)),((1),\emptyset)).$ &
\end{tabular}
\end{center}
For these $(\underline{\alpha},\underline{\beta})$, the corresponding $F_{\underline{\alpha},\underline{\beta}}$ functions can be calculated as follows.

\begin{enumerate}
\item For $\underline{\alpha}=((0),\emptyset)$ and $\underline{\beta}=((1),\emptyset)$, we obtain the map
 $$
 F_{((0),\emptyset),((1),\emptyset)}:NG_{1,1}\times NG_{1,1} \longrightarrow NG_{2,1}.
 $$
This map can be given  by
$$
F_{((0),\emptyset),((1),\emptyset)}(x, y)=[s_0^{h^{(11)}}(x),s_1^{h^{(11)}}(y)][s_1^{h^{(11)}}(y),s_1^{h^{(11)}}(x)]\in NG_{2,1}
$$
for $x,y\in NG_{1,1}$.

\item For $\underline{\alpha}=((1),\emptyset), \underline{\beta}=(\emptyset,(0))$, we get the map
$$
 F_{((1),\emptyset),(\emptyset,(0))}:NG_{1,1}\times NG_{2,0} \longrightarrow NG_{2,1}
$$
given by
$$
F_{((1),\emptyset),(\emptyset,(0))}(x, a)=[s_1^{h^{(11)}}(x),s_0^{v^{(02)}}(a)]\in NG_{2,1}
$$
for $x\in NG_{1,1}$ and $a\in NG_{2,0}$.

\item
For $\underline{\alpha}=((0),\emptyset), \underline{\beta}=(\emptyset,(0))$, we get the map
$$
F_{((0),\emptyset),(\emptyset,(0))}:NG_{1,1}\times NG_{2,0} \longrightarrow NG_{2,1}
$$
given by
$$
F_{((0),\emptyset),(\emptyset,(0))}(x, a)=[s_0^{h^{(11)}}(x),s_0^{v^{(02)}}(a)]\in NG_{2,1}
$$
for $x\in NG_{1,1}$ and $a\in NG_{2,0}$.

\item
For $\underline{\alpha}=((1),(0)), \underline{\beta}=((0),\emptyset)$, we get the following map
$$
 F_{((1),(0)),((0),\emptyset)}:NG_{1,0}\times NG_{1,1} \longrightarrow NG_{2,1}.
$$
It is given by
$$
F_{((1), (0)),((0),\emptyset)}(x, y)=[s_1^{h^{(11)}}s_0^{v^{(10)}}(x),s_0^{h^{(11)}}(y)]\in NG_{2,1}
$$
for $x\in NG_{1,0}$ and $y\in NG_{1,1}$.

\item For $\underline{\alpha}=((0),(0)), \underline{\beta}=((1),\emptyset)$, we get the following map
$$
 F_{((0),(0)),((1),\emptyset)}:NG_{1,0}\times NG_{1,1} \longrightarrow NG_{2,1}.
$$
This can be defined by
$$
F_{((0), (0)),((1),\emptyset)}(x,y)=[s_0^{h^{(11)}}s_0^{v^{(10)}}(x),s_1^{h^{(11)}}(y)]\in NG_{2,1}
$$
for $x\in NG_{1,0}$ and $y\in NG_{1,1}$.
\end{enumerate}

Let $(n,m)=(2,2)$.  By
choosing appropriate $\underline{\alpha },\underline{\beta }$ from the set $%
S(2)\times S(2)$, we can calculate the non identity maps with
codomain $NG_{2,2}$. The possible $\underline{%
\alpha },\underline{\beta }$ are given as follows.

\begin{center}
\begin{tabular}{ll}
1. $(\underline{\alpha },\underline{\beta })=(((0),\emptyset
),((1),\emptyset ))$ & 2. $(\underline{\alpha },\underline{\beta }%
)=(((1),\emptyset ),(\emptyset ,(0)))$ \\
3. $(\underline{\alpha },\underline{\beta })=(((0),\emptyset ),(\emptyset
,(0)))$ & 4. $(\underline{\alpha },\underline{\beta })=(((1),(0)),((0),%
\emptyset ))$ \\
5. $(\underline{\alpha },\underline{\beta })=(((0),(0)),((1),\emptyset )).$
& 6. $(\underline{\alpha },\underline{\beta })=((\emptyset ,(0)),(\emptyset
,(1)))$%
\end{tabular}
\end{center}

For these $(\underline{\alpha},\underline{\beta})$, the corresponding $F_{%
\underline{\alpha},\underline{\beta}}$ functions can be calculated similarly. These functions were listed in the table of Section \ref{table}.

\subsection{Appendix B}

\textbf{The Proof of Proposition \ref{2}:}
The action of $ a\in NG_{0,0}$ on $ (x,y)\in NG_{0,1}\times NG_{1,0}$ can be given
by
\[
^{a}(x,y)=(s_{0}^{v^{(00)}}(a)xs_{0}^{v^{(00)}}(a)^{-1},s_{0}^{h^{(00)}}(a)ys_{0}^{h^{(00)}}(a)^{-1}).
\]

$CM1)$ For $ a\in NG_{0,0} $ and $ x\in NG_{0,1}, y\in
NG_{1,0},$ we obtain
\begin{eqnarray*}
\partial (^{a}(x,y)) & =&\partial (s_{0}^{v^{(00)}}(a)xs_{0}^{v^{(00)}}(a)^{-1},s_{0}^{h^{(00)}}(a)ys_{0}^{h^{(00)}}(a)^{-1}) \\
&=&d_{1}^{v^{(01)}}(s_{0}^{v^{(00)}}(a)xs_{0}^{v^{(00)}}(a)^{-1}) d_{1}^{h^{(1,0)}}(s_{0}^{h^{(00)}}(a)ys_{0}^{h^{(00)}}(a)^{-1}) \\
& =&(ad_{1}^{v^{(01)}}(x)a^{-1})(ad_{1}^{h^{(10)}}(y)a^{-1}) \\
& =&ad_{1}^{v^{(01)}}(x)d_{1}^{h^{(10)}}(y)a^{-1} \\
& =&a\partial (x,y)a^{-1}.
\end{eqnarray*}

$CM2)$ For $(x_{1},y_{1}),(x_{2},y_{2})\in NG_{0,1}\times NG_{1,0},$
\begin{eqnarray*}
^{\partial (x_{1},y_{1})}(x_{2},y_{2}) & =&
^{{d_{1}^{{^{v}}^{(0,1)}}}
(x_{1}){d_{1}^{h^{(1,0)}}(y_{1})}}(x_{2},y_{2}) \\
& =&\bigg(s_{0}^{{^{v}}^{(00)}}\big(d_{1}^{{^{v}}^{(01)}}
(x_{1})d_{1}^{h^{(10)}}(y_{1})\big)x_{2}s_{0}^{{^{v}}^{(00)}}\big(
d_{1}^{{^{v}}^{(01)}}(x_{1})d_{1}^{h^{(10)}}(y_{1})\big)^{-1}, \\
&&
s_{0}^{h^{(00)}}\big(d_{1}^{{^{v}}^{(01)}}(x_{1})d_{1}^{h^{(10)}}(y_{1})
\big)y_{2}s_{0}^{h^{(00)}}\big(d_{1}^{{^{v}}^{(01)}}(x_{1})d_{1}^{h^{(10)
}}(y_{1})\big)^{-1}\bigg) \\
&
=&\bigg(s_{0}^{{^{v}}^{(00)}}d_{1}^{{^{v}}^{(01)}}(x_{1})
d_{1}^{h^{(11)}}s_{0}^{v^{(01)}}(y_{1})x_{2}\big(d_{1}^{h^{(11)}}s_{0}^{v^{(10)}}y_{1}^{-1}\big)
s_{0}^{v^{(00)}}d_{1}^{v^{(01)}}(x_{1})^{-1},\\
&&\big(d_{1}^{v^{(11)}}s_{0}^{h^{(01)}}x_{1}\big)s_{0}^{h^{(00)}}
d_{1}^{h^{(10)}}(y_{1})y_{2}\big(s_{0}^{h^{(00)}}d_{1}^{h^{(10)}}y_{1}^{-1}
\big)\big(d_{1}^{v^{(11)}}s_{0}^{h^{(01)}}x_{1}^{-1}\big)\bigg)\\
&=&\bigg(s_{0}^{{^{v}}^{(00)}}d_{1}^{{^{v}}^{(01)}}(x_{1})x_{2}s_{0}^{v^{(00)}}d_{1}^{v^{(01)}}(x_{1})^{-1}, s_{0}^{h^{(00)}}d_{1}^{h^{(10)}}(y_{1})y_{2} s_{0}^{h^{(00)}}d_{1}^{h^{(10)}}y_{1}^{-1}\bigg).
\end{eqnarray*}
Since $NG_{1,1}=\{1\}$ for $x_1, x_{2}\in NG_{0,1}$  we have
\begin{align*}
\partial_2^{v^{(02)}}(F_{(\emptyset, (0)),(\emptyset, (1))}(x_1,x_2)&=[s_{0}^{v^{(00)}}d_{1}^{v^{(01)}}(x_{1}),x_2][x_2,x_1]\\
&=s_{0}^{v^{(00)}}d_{1}^{v^{(01)}}(x_{1})x_2s_{0}^{v^{(00)}}d_{1}^{v^{(01)}}(x_{1})^{-1}x_1 x_2^{-1}x_1^{-1}\in \partial_2^{v^{(02)}}(NG_{0,2}\cap D_{0,2})=\{1\},\\
\end{align*}
and we obtain
$$
s_{0}^{v^{(00)}}d_{1}^{v^{(01)}}(x_{1})x_2s_{0}^{v^{(00)}}d_{1}^{v^{(01)}}(x_{1})^{-1}=x_1 x_2x_1^{-1}.
$$
Similarly, for $y_1, y_{2}\in NG_{1,0}$  we have
\begin{align*}
\partial_2^{h^{(20)}}(F_{((0),\emptyset),((1),\emptyset)}(y_1,y_2)&=[s_{0}^{h^{(00)}}d_{1}^{h^{(10)}}(y_{1}),y_2][y_2,y_1]\\
&=s_{0}^{h^{(00)}}d_{1}^{h^{(10)}}(y_{1})y_2s_{0}^{h^{(00)}}d_{1}^{h^{(10)}}(y_{1})^{-1}y_1 y_2^{-1}y_1^{-1}\in \partial_2^{h^{(20)}}(NG_{2,0}\cap D_{2,0})=\{1\},\\
\end{align*}
and we obtain
$$
s_{0}^{h^{(00)}}d_{1}^{h^{(10)}}(y_{1})y_2s_{0}^{h^{(00)}}d_{1}^{h^{(10)}}(y_{1})^{-1}=y_1 y_2 y_1^{-1}.
$$
Thus we get
\begin{align*}
^{\partial(x_{1},y_{1})}(x_{2},y_{2})&=(x_1 x_2x_1^{-1},y_1 y_2 y_1^{-1})\\
&=(x_1,y_1)(x_2,y_2)(x_1,y_1)^{-1}
\end{align*}
and this is the second condition of crossed module. $\Box$

\textbf{The Proof of Proposition \ref{3}:}
The actions
of $ NG_{0,1} $ and $ NG_{1,0} $ on $ NG_{1,1} $ are given by
\begin{enumerate}
\item $^{y}x=s_{0}^{{v }^{(00)}}d_{1}^{{h }^{(10)}}(y)xs_{0}^{{v
}^{(00)}}d_{1}^{{h }^{(10)}}(y^{-1})
=d_{1}^{h^{(11)}}s_{0}^{v^{(10)}}(y)(x)d_{1}^{h^{(11)}}s_{0}^{v^{(10)
}}(y)^{-1}$
\item $^{x}y=s_{0}^{{h }^{(00)}}d_{1}^{{v }^{(01)}}(x)(y)s_{0}^{{h}^{(00)}}d_{1}^{{v}^{(01)}}({x}^{-1})=d_{1}^{v^{(11)}}s_{0}^{h^{(10)}}(x)(y)d_{1}^{v^{(11)}}s_{0}^{h^{(10)}}(x)^{-1}$
\item $^{y}z=s_{0}^{v^{(10)}}(y)(z)s_{0}^{v^{(10)}}(y)^{-1} $
\item $ ^{x}z=s_{0}^{h^{(01)}}(x)(z)s_{0}^{h^{(01)}}(x)^{-1}$
\end{enumerate}
for $x\in NG_{0,1},\ y\in NG_{1,0} $ and $ z\in NG_{1,1}.$

$(i)$ Since $NG_{1,2}=NG_{2,1}=NG_{0,2}=NG_{2,0}=\{1\}$, by using the method given in the previous proposition, it can be easily shown that the maps $\partial_1^{h^{(10)}},\partial_1^{v^{(01)}},\partial_1^{h^{(11)}}$ and $\partial_1^{v^{(11)}}$ are crossed modules.

$(iv)$ For $x\in N_{0,1} $ and $y\in NE_{1,0}$, we obtain
\begin{align*}
\partial_1^{h^{(11)}}h(x,y)&=\partial_1^{h^{(11)}}(s_0^{h^{(01)}}(x)s_0^{v^{(10)}}(y)s_0^{h^{(01)}}(x)^{-1})s_0^{v^{(10)}}(y)^{-1}\\
&=x\partial_1^{h^{(11)}}s_0^{v^{(10)}}(y)x^{-1}\partial_1^{h^{(11)}}s_0^{v^{(10)}}(y)^{-1}\\
&=xs_0^{v^{(00)}}d_1^{h^{(10)}}(y)x^{-1}s_0^{v^{(00)}}d_1^{h^{(10)}}(y)^{-1}\\
&=x ^yx^{-1}.
\end{align*}

$(v)$ We obtain for $x\in NE_{0,1}$ and $y\in NE_{1,0}$
\begin{align*}
\partial_1^{v^{(11)}}h(x,y)&=\partial_1^{v^{(11)}}[s_0^{h^{(01)}}(x),s_0^{v^{(10)}}(y)]\\
&=[\partial_1^{v^{(11)}}s_0^{h^{(01)}}(x),y]\\
&=[s_0^{h^{(00)}}d_1^{v^{(01)}}(x),y]\\
&=s_0^{h^{(00)}}d_1^{v^{(01)}}(x)ys_0^{h^{(00)}}d_1^{v^{(01)}}(x)^{-1}y^{-1}\\
&=^x y y^{-1}.
\end{align*}

$(vi)$ For $z\in NG_{1,1}$ and $y\in NG_{1,0}$, we get
$$
h(\partial_{1}^{h^{(11)}}(z),y) =[s_{0}^{h^{(01)}}d_{1}^{h^{(11)}}(z),s_{0}^{v^{(10)}}(y)].\\
$$
On the other hand for $z\in NG_{1,1}$ and $s_0^{v^{(10)}}y\in NG_{1,1}$, we obtain also
$$
F_{((0),\emptyset),((1),\emptyset)}(z,s_0^{v^{(10)}}y)=[s_{0}^{h^{(11)}}z,s_{1}^{h^{(11)}}s_{0}^{v^{(10)}}(y)][s_{1}^{h^{(11)}}s_{0}^{v^{(10)}}(y), s_{1}^{h^{(11)}}z]\in NG_{2,1}=\{1\},
$$
and
\begin{align*}
\partial_2^{h^{(21)}}F_{((0),\emptyset),((1),\emptyset)}(z,s_0^{v^{(10)}}y)&=[s_{0}^{h^{(01)}}d_{1}^{h^{(11)}}(z),s_{0}^{v^{(10)}}(y)]
[s_{0}^{v^{(10)}}(y), z]=1.
\end{align*}
Thus we get
\begin{align*}
h(\partial_{1}^{h^{(11)}}(z),y)& =[s_{0}^{h^{(01)}}d_{1}^{h^{(11)}}(z),s_{0}^{v^{(10)}}(y)]\\
&=[z,s_{0}^{v^{(10)}}(y)]\\
&=zs_{0}^{v^{(10)}}(y)z^{-1}s_{0}^{v^{(10)}}(y)^{-1}\\
&=z ^y(z^{-1}).
\end{align*}

$(vii)$ For $z\in NG_{1,1}$ and $y\in NG_{0,1}$, we get
$$
h(y,\partial_{1}^{v^{(11)}}(z)) =[s_{0}^{h^{(01)}}(y),s_{0}^{v^{(10)}}d_{1}^{v^{(11)}}(z)].\\
$$
On the other hand for $z\in NG_{1,1}$ and $s_0^{h^{(01)}}y \in NG_{1,1}$, we obtain also
$$
F_{(\emptyset,(0)),(\emptyset,(1))}(z,s_0^{h^{(01)}}y)=[s_{0}^{v^{(11)}}z,s_{1}^{v^{(11)}}s_{0}^{h^{(01)}}(y)][s_{1}^{v^{(11)}}s_{0}^{h^{(01)}}(y), s_{1}^{v^{(11)}}z]\in NG_{1,2}=\{1\},
$$
and
\begin{align*}
\partial_2^{v^{(12)}}F_{(\emptyset,(0)),(\emptyset,(1))}(z,s_0^{h^{(01)}}y)&=[s_{0}^{v^{(10)}}d_{1}^{v^{(11)}}(z),s_{0}^{h^{(01)}}(y)]
[s_{0}^{h^{(01)}}(y), z]=1.
\end{align*}
Thus we get
\begin{align*}
h(y,\partial_{1}^{v^{(11)}}(z))& =[s_{0}^{h^{(01)}}y, s_{0}^{v^{(10)}}d_{1}^{v^{(11)}}(z)]\\
&=[s_{0}^{h^{(01)}}(y),z]\\
&=s_{0}^{h^{(01)}}(y)zs_{0}^{h^{(01)}}(y)^{-1}z^{-1}\\
&= ^y(z)z^{-1}.
\end{align*}
$\Box$

\textbf{The Proof of Proposition \ref{5}}:
In each direction, by using Proposition \ref{4} and the images of the functions  $F_{\underline{\alpha},\underline{\beta}}$ similarly to \cite{mutpor2}, the proof can be easily given. For example, one has:

1.\
\begin{eqnarray*}
\partial _{2}^{h^{(2q)}}\{x,y\}
&=&d_{2}^{h^{(2q)}}[s_{1}^{h^{(1q)}}(x),s_{1}^{h^{(1q)}}(y)][s_{1}^{h^{(1q)}}(y),s_{0}^{h^{(1q)}}(x)]
\\
&=&[x,y][y,d_{2}^{h^{(2q)}}s_{0}^{h^{(1q)}}(x)]
\\
&=&xyx^{-1}s_{0}^{h^{(0q)}}d_{1}^{h^{(1q)}}(x)(y)^{-1}s_{0}^{h^{(0q)}}d_{1}^{h^{(1q)}}(x^{-1})
\\
&=&xyx^{-1}(^{\partial _{1}^{h^{(1q)}}(x)}y)^{-1}.
\end{eqnarray*}

2.\

Since
\begin{align*}
\partial_3^{h^{(3q)}}\left(F_{((0),\emptyset),((1),\emptyset)}(z,z')\right)^{-1}
&=[z',z][s_{1}^{h^{(1q)}}d_{2}^{h^{(2q)}}(z),s_{1}^{h^{(1q)}}d_{2}^{h^{(2q)}}(z')]\\
&\ \ \hspace{2.cm}[s_{1}^{h^{(1q)}}d_{2}^{h^{(2q)}}(z'),s_{0}^{h^{(1q)}}d_{2}^{h^{(2q)}}(z)]
\in\partial_3^{h^{(3q)}}(NG_{3,q}\cup D_{3,q})=\{1\} \\
\end{align*}
 for $z,z'\in NG_{2,q}$, we have
\begin{align*}
\{\partial _{2}^{h^{(2q)}}z,\partial _{2}^{h^{(2q)}}z'\}
&=[s_{1}^{h^{(1q)}}d_{2}^{h^{(2q)}}(z),s_{1}^{h^{(1q)}}d_{2}^{h^{(2q)}}(z')][s_{1}^{h^{(1q)}}d_{2}^{h^{(2q)}}(z'),s_{0}^{h^{(1q)}}d_{2}^{h^{(2q)}}(z)]\\
&=[z,z'].
\end{align*}
$\Box$

\bigskip
$
\begin{array}{llllllll}
 &     & \text{\"O. G{\"u}rmen Alansal and E. Ulualan} &  \\
  &  & \text{Dumlup\i nar University} &  \\
 &   & \text{Science and Art Faculty%
} &  \\
   &  & \text{%
Mathematics Department} &  \\
 &  & \text{K\"{u}tahya,
TURKEY} &  \\
  &  & \text{ogurmen@gmail.com}&\\
&&\text{eulualan@gmail.com} &
\end{array}%
$
\end{document}